\documentclass[11pt]{amsart2000}
\newtheorem{theorem}{Theorem}
\newtheorem{corollary}{Corollary}
\newtheorem*{notheorem}{Proposition}
\theoremstyle{definition}
\newtheorem{definition}{Defintion}
\newtheorem*{nodefinition}{Definition}
\newtheorem*{remark}{Remark}

\begin{document}
\setlength{\unitlength}{0.01in}
\linethickness{0.01in}
\begin{center}
\begin{picture}(474,66)(0,0)
\multiput(0,66)(1,0){40}{\line(0,-1){24}}
\multiput(43,65)(1,-1){24}{\line(0,-1){40}}
\multiput(1,39)(1,-1){40}{\line(1,0){24}}
\multiput(70,2)(1,1){24}{\line(0,1){40}}
\multiput(72,0)(1,1){24}{\line(1,0){40}}
\multiput(97,66)(1,0){40}{\line(0,-1){40}}
\put(143,66){\makebox(0,0)[tl]{\footnotesize Proceedings of the Ninth Prague Topological Symposium}}
\put(143,50){\makebox(0,0)[tl]{\footnotesize Contributed papers from the symposium held in}}
\put(143,34){\makebox(0,0)[tl]{\footnotesize Prague, Czech Republic, August 19--25, 2001}}
\end{picture}
\end{center}
\vspace{0.25in}
\setcounter{page}{205}
\title[Proper and admissible topologies in closure spaces]{Proper and
admissible topologies in the setting of closure spaces}
\author{Mila Mr\v{s}evi\'c}
\address{Faculty of Mathematics, University of Belgrade\\
Beograd, Yugoslavia}
\email{mrsevic@matf.bg.ac.yu}
\thanks{This article will be expanded and submitted for publication
elsewhere.}
\thanks{Mila Mr\v{s}evi\'c,
{\em Proper and admissible topologies in the setting of closure spaces},
Proceedings of the Ninth Prague Topological Symposium, (Prague, 2001),
pp.~205--216, Topology Atlas, Toronto, 2002}
\begin{abstract}
A \v Cech closure space $(X,u)$ is a set $X$ with a (\v Cech) closure
operator $u$ which need not be idempotent. 
Many properties which hold in topological spaces hold in \v Cech closure
spaces as well.

The notions of proper (splitting) and admissible (jointly 
continuous) topologies are introduced on the sets of continuous functions
between \v Cech closure spaces. 
It is shown that some well-known results of Arens and Dugundji \cite{AD}
and Iliadis and Papadopoulos \cite{IP} are true in this setting.

We emphasize that Theorems 1--10 encompass the results of A. di Concilio
\cite{C} and Georgiou and Papadopoulos \cite{GP1, GP2} for the spaces
of continuous-like functions as $\theta$-continuous, strongly and weakly
$\theta$-contin\-uous, weakly and super-continuous.
\end{abstract}
\maketitle

\section{\v Cech closure spaces}

An operator $u : \mathcal{P}(X) \rightarrow \mathcal{P}(X)$ defined on the 
power set $\mathcal{P}(X)$ of a set $X$ satisfying the axioms:
\begin{itemize}
\item[(C1)]
$u(\emptyset) = \emptyset$,
\item[(C2)]
$A \subset u(A)$ for every $A \subset X$,
\item[(C3)]
$u(A \cup B) = u(A) \cup u(B)$ for all $A,B \subset X$,
\end{itemize}
is called a {\it \v Cech closure operator} and the pair $(X,u)$ is a
{\it \v Cech closure space}. 
For short, the space will be noted by $X$ as well, and called a 
{\it closure space}.

A subset $A$ is {\it closed} in the closure space $(X,u)$ if $u(A)=A$ holds. 
It is {\it open} if its complement is closed. 
The empty set and the whole space are both open and closed.

The \space {\it interior operator}
${\rm int}_u:\mathcal{P}(X) \rightarrow \mathcal{P}(X)$ 
is defined by means of the closure operator in the usual way: 
${\rm int}_u ={\rm c}\circ u\circ {\rm c}$, where 
${\rm c}:\mathcal{P}(X) \rightarrow \mathcal{P}(X)$
is the complement operator. 
A subset $U$ is a {\it neighbourhood} of a point $x$ (subset $A$) in $X$
if $x\in {\rm int}_uU (A\subset {\rm int}_uU )$ 
holds. 
We denote by $\mathcal{N}(x)$ the collection of all neighbourhoods 
(the {\it neighbourhood system)} at the point $x$. 
By (C3), the intersection of two (and thus finitely many) neighbourhoods
at $x$ is a neighbourhood at $x$ again. 
The condition (C1) is equivalent to ${\rm int}_uX = X$, that is to
$X\in\mathcal{N}(x)$ for every $x\in X$, and ${\rm int}_u A \subset A$
for every $A \subset X$ is equivalent to (C2).

In a closure space $(X,u)$ a family 
$\mathcal{U}(x) \subset \mathcal{N}(x)$ 
is a neighbourhood (local) base at a point $x$ if the following
axioms are satisfied:
\begin{itemize}
\item[(Nb1)]
$\mathcal{U}(x) \ne \emptyset$ for every $x \in X$,
\item[(Nb2)]
$x \in U$ for every $U \in \mathcal{U}(x)$,
\item[(Nb3)]
$U_{1}, U_{2} \in \mathcal{U}(x) \Rightarrow 
(\exists U \in \mathcal{U}(x)) U \subset U_{1} \cap U_{2}$. 
\end{itemize}

A family $\mathcal{U}(x) \subset \mathcal{N}(x)$
is a neighbourhood (local) subbase at a point
$x$ if the conditions (Nb1) and (Nb2) are fulfilled.

If a collection $\{\mathcal{U}(x) \mid x \in X\}$ of filters on $X$
satisfies the conditions (Nb1)--(Nb3), then there is exactly one closure
operator $u$ for $X$ such that $\mathcal{U}(x)$ is a neighbourhood base at
$x$ for each $x\in X$. 
The operator $u$ is defined by: 
$$
u(A) = 
\{ x\in X \mid U\in \mathcal{U}(x) \Rightarrow U \cap A \ne \emptyset \}.
$$

Let $(X,u_{1})$ and $(X,u_{2})$ be closure spaces. 
The closure $u_{1}$ is {\it coarser} than the closure $u_{2}$, or $u_{2}$
is {\it finer} than $u_{1}$, denoted by $u_{1} \leq u_{2}$, if 
$u_{1}{(A) \supset u}_{2}(A)$ for every $A \subset X$. 
So defined relation $\leq$ is a partial order on the set of all closure
spaces.

Let $\{ u_{\alpha }\}$ be a collection of closure operators on a set $X$. 
The infimum (meet) and supremum (join) operators for $\{ u_{\alpha }\}$
are the operators
$u_{0} = {\land u}_{\alpha }$ and
$u = {\lor u}_{\alpha }$ respectively, defined by:
$\mathcal{U}_{0}(x) = \bigcap_{\alpha} \mathcal{N}_{\alpha}(x)$ 
is a neighbourhood base (system) and 
$\mathcal{U}(x) = \bigcup_{\alpha} \mathcal{N}_{\alpha}(x)$ 
is a neighbourhood subbase at $x\in X$, for $u_{0}$ and $u$ respectively.

Many topological notions can be defined in the class of closure spaces by
means of neighbourhoods.

Let ${\rm M} $ be a directed set and $(x_{\mu})_{\mu \in {\rm M}}$ a net
in $(X,u)$.
The net $(x_{\mu})$ converges to a point $x \in X$ if for every
neighbourhood $U$ of $x$ there is a $\mu \in {\rm M}$ such that for every
$\mu' \in {\rm M}$, $\mu' \geq \mu \Rightarrow x_{\mu '} \in U$.
Similarly, $x$ is {\it an accumulation point} of the net $(x_{\mu})$ if
for every neighbourhood $U$ of $x$ and every $\mu \in {\rm M}$ there is a
$\mu' \in {\rm M}$ such that $\mu' \geq \mu$ and $x_{\mu'} \in U$. 
For every point $x$ the neighbourhood system $\mathcal{N}(x)$ is a filter
on $X$ such that $x\in \bigcap \mathcal{N}(x)$.
Moreover it is a set directed by the inverse inclusion $\supset$ and
every net $(x_{U})_{U\in \mathcal{N}(x)}$ with $x_{U} \in U$,
converges to $x$.

Let $(X,u)$ and $(Y,v)$ be two closure spaces.
A function $f:(X,u)\rightarrow (Y,v)$ is {\it continuous at} $x \in X$ if
``close points are mapped into close ones'', that is if the following 
holds
$$
A \subset X \land x \in u(A) \Rightarrow f(x) \in v(f(A)).
$$
This condition is equivalent to:
\begin{itemize}
\item[(i)]
the inverse image of every neighbourhood of $f(x)$ is a neighbourhood of 
$x$;
\item[(ii)]
for every net $(x_{\mu})$ that converges to $x$, the net $(f(x_{\mu}))$
converges to $f(x)$;
\item[(iii)]
if $x$ is an accumulation point of a net $(x_{\mu}), f(x)$ is an
accumulation point of the net $(f(x_{\mu}))$.
\end{itemize}

A function $f:(X,u) \rightarrow (Y,v)$ is {\it continuous} if it is
continuous at every point of $X$. 
This condition is equivalent to: 
\begin{itemize}
\item[(i)]
$f(u(A)) \subset v(f(A))$ for every $A \subset X$;
\item[(ii)]
$u(f^{-1}(B)) \subset f^{-1}(v(B))$ for every $B \subset Y$.
\end{itemize}

The {\it product of a family 
$\{ (X_{\alpha}, u_{\alpha}) \mid \alpha \in {\rm A} \}$
of closure spaces}, denoted by $\Pi (X_{\alpha} ,u_{\alpha})$, 
is the set 
$X = \Pi_{\alpha \in {\rm A}}X_{\alpha}$
endowed with the closure operator $u$ defined by means of neighbourhoods:
for every $x \in X$ the family 
$$\mathcal{U}(x) = 
\{ \pi_{\alpha}^{-1} (V) \mid 
\alpha \in {\rm A}, V \in \mathcal{N}_{\alpha}(x_{\alpha }) \}$$
is a neighbourhood subbase at $x$ in $(X,u)$.
Here $\pi_{\alpha }$ are the projections, while 
$\mathcal{N}_{\alpha}(x_{\alpha })$ is the neighbourhood system at 
$x_{\alpha} = \pi_{\alpha}(x)$ in $X_{\alpha}$.
There exists exactly one closure operator $u$ such that $\mathcal{U}(x)$
is a local subbase at $x$ in $(X,u)$ for every $x \in X$. 
Canonical neighbourhoods of $x$ are of the form 
$\bigcap_{i=1}^{k} \pi_{\alpha _{i}}^{-1}(V_{i})$.

The projections are continuous as well as the restrictions of continuous
functions. 
The composition of two continuous mappings is continuous and a mapping 
$f:(X,u) \rightarrow \Pi (Y_{\alpha}, v_{\alpha})$ 
is continuous at $x \in X$ if and only if each composition 
$\pi_{\alpha}\circ f$ is continuous at $x$. 
Also the product 
$f: \Pi (X_{\alpha} ,u_{\alpha}) \rightarrow \Pi (Y_{\alpha}, v_{\alpha})$
of a family $\{ f_{\alpha }\}$ of continuous mappings is continuous.

A well-known example of a \v Cech closure operator which is not a 
Kuratowski closure operator in general, is the so called $\theta$-closure.
It was defined by Veli\v cko \cite{V} in the following way:
Let $(X,\mathcal{T})$ be a topological space and let $A \subset X$. 
A point $x \in X$ is in the $\theta$-{\it closure} of $A$, denoted by
${\rm cl}_{\theta}A$ (or $\mathcal{T}{\rm cl}_{\theta}A$), if each
closed neighbourhood of $x$ intersects $A$. 
Neighbourhood bases in $(X,{\rm cl}_{\theta})$ consist of closed 
neighbourhoods (or closures of open neighbourhoods) in $(X,\mathcal{T})$
at every point $x$.

Let $(X,\mathcal{T})$ be the product space of a family 
$\{(X_{\alpha},\mathcal{T}_{\alpha})\}$ of topological spaces. 
The $\theta$-closure space of $(X,\mathcal{T})$ is the product of the
$\theta$-closure spaces of $(X_{\alpha},\mathcal{T}_{\alpha})$, 
i.e.\ 
$(X,\mathcal{T}{\rm cl}_{\theta}) = 
\Pi (X_{\alpha},\mathcal{T}_{\alpha}{\rm cl}_{\theta})$.

A function $f:(X,\mathcal{T}) \rightarrow (Y,\mathcal{V})$ is
$\theta$-{\it continuous at} $x \in X$ if for every neighbourhood $V$ of
$f(x)$ there is a neighbourhood $U$ of $x$ such that 
$f(\overline{U}) \subset \overline{V}$.
A function $f:(X,\mathcal{T}) \rightarrow (Y,\mathcal{V})$ is 
$\theta$-{\it continuous} if it is $\theta$-continuous at each of its
points. 
Every continuous function is $\theta$-continuous, but the converse does
not hold in general.

$\theta$-continuity is not a continuity concept in the class of 
topological spaces, but it is in the class of \v Cech closure spaces. 
Namely,
\begin{notheorem}
A function $f:(X,\mathcal{T}) \rightarrow (Y,\mathcal{V})$
is $\theta$-continuous if and only if the function 
$f:(X,\mathcal{T}{\rm cl}_{\theta}) \rightarrow 
(Y,\mathcal{V}{\rm cl}_{\theta})$ is a
continuous mapping of (\v Cech) closure spaces.
\end{notheorem}

Hence the following characterizations of $\theta$-continuity:
\begin{notheorem}
A function $f:(X,\mathcal{T})\rightarrow (Y,\mathcal{V})$
is $\theta$-continuous if and only if:
\begin{enumerate}
\item
$f(\mathcal{T}{\rm cl}_{\theta}A) \subset 
\mathcal{V}{\rm cl}_{\theta}f(A)$ 
for every $A\subset X$;
\item
$\mathcal{T}{\rm cl}_{\theta} f^{-1}(B) \subset 
f^{-1}(\mathcal{V}{\rm cl}_{\theta}B)$ 
for every $B\subset Y$.
\end{enumerate}
\end{notheorem}

The next statement follows from the definitions and the properties of
$\theta$-closure.
\begin{notheorem}
Let $X, Y, Z$ be topological spaces.
A function $g:Z\times X\rightarrow Y$ is $\theta$-continuous if and only
if the function 
$g:(Z,{\rm cl}_{\theta}) \times (X,{\rm cl}_{\theta}) \rightarrow 
(Y,{\rm cl}_{\theta})$ is continuous.
\end{notheorem}

A closure space $(X,u)$ is:
\begin{itemize}
\item[(i)]
{\it regular} if for each point $x$ and each subset $A$ such that 
$x \notin u(A)$, there exist neighbourhoods $U$ of $x$ and $V$ of $A$ such
that $U \cap V = \emptyset$;
\item[(ii)]
{\it compact} \space if each net in $(X,u)$ has an accumulation point.
\end{itemize}

A closure space $(X,u)$ is regular if and only if for each point $x$ and
each neighbourhood $U$ of $x$, there is a neighbourhood $U_{1}$ of $x$
such that $u(U_{1})\subset U$.

Compactness can be characterized by means of covers.
\cite[41 A.9. Theorem]{vC}
An {\it interior cover} of $(X,u)$ is a cover $\{ G_{\alpha }\}$
such that the collection $\{{\rm int}_{u}G_{\alpha }\}$ covers $X$. 
The space is compact if and only if every interior cover has a
finite subcover.

We give the following

\begin{nodefinition}
A collection $\{ G_{\alpha }\}$ is an {\it interior cover} of a set
$A$ in $(X,u)$ if the collection $\{ {\rm int}_{u}G_{\alpha }\}$ covers
$A$. 
A subset $A$ is {\it compact} if every interior cover of $A$ has a finite
subcover.
\end{nodefinition}

All notions not explained here can be found in \cite{vC}.

\section{Proper and admissible topologies in the setting of closure
spaces}

Let $X,Y$ and $Z$ be three nonempty sets.
For every function $g:Z\times X \rightarrow Y$ there is a function 
${\rm E}(g)$ or $g^{\ast}$ from $Z$ to $Y^X$, the set of all functions
from $X$ to $Y$, defined by $(g^{\ast}(z))(x)=g(z,x)$. 
The mapping ${\rm E}:Y^{Z\times X} \rightarrow (Y^{X})^{Z}$ is called the
{\it exponential function}. 
By $\varepsilon $ we denote the {\it evaluation mapping} from 
${Y}^{X}\times X$ to $Y$ defined by $\varepsilon(f,x) = f(x)$.

If $X,Y$ and $Z$ are topological or closure spaces, in particular sets of
continuous functions can be considered. 
Now on $Y^X$ will mean the set of {\it all continuous functions} from 
$X$ to $Y$. 
The set $Y^X$ can be endowed with different topologies. 
The question is: 
Find the topologies on the set of functions such that 
\begin{enumerate}
\item
${\rm E}(g) = g^{\ast} \in (Y^{X})^{Z}$ 
for every $g \in Y^{Z\times X}$,
that is, for every continuous $g:Z\times X \rightarrow Y$ the function 
$g^{\ast}$ is continuous; 
and conversely,
\item
$g \in Y^{Z\times X}$ for every $g^{\ast} \in (Y^{X})^{Z}$, that is,
for every continuous $g^{\ast}:Z \rightarrow Y^{X}$ the function $g$ is
continuous.
\end{enumerate}

Following the definitions and notations used by Arens and Dugundji 
\cite{AD}, Kuratowski \cite{K} and Iliadis and Papadopoulos \cite{IP} for
the sets of continuous functions defined in the setting of topological
spaces, we give the following definitions.

\begin{definition}
Let $(X,u)$ be a closure space and $(A_{\lambda })_{\lambda \in \Lambda}$
be a net in $\mathcal{P}(X)$. 
The {\it upper limit} of the net $(A_{\lambda})$, denoted by 
$\overline{\lim\limits_{\Lambda }}A_{\lambda}$, is the set of all points
$x \in X$ such that for every $\lambda_{0}\in \Lambda$ and every 
neighbourhood $U$ of $x$ in $X$, there is a $\lambda \in \Lambda$ such 
that $\lambda{\geq \lambda}_{0}$ and 
$A_{\lambda} \cap U \ne \emptyset$. 
(See, for example, \cite{AD} and \cite{IP}.)
\end{definition}

\begin{definition}
Let $(X,u)$ and $(Y,v)$ be closure spaces and $Y^X$ be the collection of 
all continuous functions $f:(X,u) \rightarrow (Y,v)$. 
A closure operator $\sigma$ on $Y^X$ is called {\it proper (splitting)} if
for any closure space $(Z,w)$
\begin{itemize}
\item[(1)]
$g:(Z,w)\times (X,u) \rightarrow (Y,v)$ is continuous 
$\Rightarrow 
E(g) = g^{\ast}:(Z,w) \rightarrow (Y^{X},\sigma)$ is continuous;
\end{itemize}
$\sigma$ is called {\it admissible (jointly continuous)} if for every
space $(Z,w)$
\begin{itemize}
\item[(2)]
$g^{\ast}:(Z,w) \rightarrow (Y^{X},\sigma)$ is continuous 
$\Rightarrow 
g:(Z,w)\times (X,u) \rightarrow (Y,v)$ is continuous.
\end{itemize}
Let $f, f_{\lambda} \in Y^{X}$, $\lambda \in \Lambda$, where $\Lambda$ is
a directed set. 
The net $(f_{\lambda})$ {\it converges continuously} to $f$, denoted
by $f_{\lambda}\buildrel cc\over {\longrightarrow f}$, if
\begin{itemize}
\item[(3)]
the net 
$f_{\lambda }(x_{\mu })$, $(\lambda,\mu) \in \Lambda \times {\rm M}$,
converges to $f(x)$ in $(Y,v)$ whenever the net $(x_{\mu})$ converges to
$x$ in $(X,u)$.
\end{itemize}
The convergence of a net $(f_{\lambda})$ to $f$ in the space 
$(Y^{X},\sigma)$ will be denoted by 
$f_{\lambda }\buildrel \sigma\over {\longrightarrow f}$.
\end{definition}

The next results follow from definitions and the proofs are analogous to
the corresponding for the topological case. (See \cite{AD} and \cite{IP}).

\begin{theorem}
A closure operator $\sigma$ on $Y^{X}$ is admissible if and only if the
evaluation mapping 
$\varepsilon:(Y^{X},\sigma)\times (X,u) \rightarrow (Y,v)$ is
continuous.
\end{theorem}

\begin{theorem}
Let $(f_{\lambda})_{\lambda\in \Lambda}$ be a net in $Y^{X}$. 
The net $(f_{\lambda})_{\lambda \in \Lambda}$ converges continuously to
$f \in Y^{X}$ if and only if for every $x\in X$ and every neighbourhood
$V$ of $f(x)$ there is a neighbourhood $U$ of $x$ and a 
$\lambda _{0} \in \Lambda $ such that $f_{\lambda }(U) \subset V$ for all
$\lambda{\geq \lambda}_{0}$.
\end{theorem}

\begin{theorem}
\it Let $(f_{\lambda })_{\lambda \in \Lambda }$ be a net in $Y^{X}$. 
The net $(f_{\lambda})_{\lambda \in \Lambda }$ converges continuously to
$f\in Y^{X}$ if and only if the following holds:
\begin{itemize}
\item[(4)]
$\overline{\lim\limits_{\Lambda }} f_{\lambda}^{-1}(B) \subset f^{-1}(v(B))$ 
for every subset $B$ in $Y$.
\end{itemize}
\end{theorem}

\begin{proof}
Let $f_{\lambda }\buildrel cc\over {\longrightarrow f}$, $B \subset Y$
and $x \in \overline{\lim\limits_{\Lambda }} f_{\lambda}^{-1}(B)$. 
For every neighbourhood $V$ of $f(x)$, by Theorem 2, there is a 
neighbourhood $U$ of $x$ and a $\lambda_{0}\in \Lambda$ such that 
$f_{\lambda}(U)\subset V$ for all $\lambda{\geq \lambda}_{0}$. 
By the assumption, there is a $\lambda' \geq \lambda_{0}$ such that
$f_{\lambda'}^{-1}(B) \cap U \ne \emptyset$. 
It follows that $f_{\lambda'}(U) \cap B \ne \emptyset$. 
Thus 
$$(\forall V \in \mathcal{N}(f(x)))\ V \cap B \ne \emptyset$$
implies $f(x) \in v(B)$, 
hence $x \in f^{-1}(v(B))$.

For the converse, let $(f_{\lambda})_{\lambda \in \Lambda }$ be a net in
$Y^{X}$ such that the condition (4) holds. 
For every $x \in X$ and every neighbourhood $V$ of $f(x)$, 
$$f(x) \in {\rm int}_vV = {\rm c}(v({\rm c}(V)))
\mbox{ implies } 
f(x) \notin v(V^{\rm c}).$$
Thus 
$x \notin f^{-1}(v(V^{\rm c}))$ 
implies
$x \notin \overline{\lim\limits_{\Lambda}} f_{\lambda}^{-1}(V^{\rm c})$. 
Hence 
$$\begin{array}{lll}
(\exists U \in \mathcal{N}(x)) 
(\exists\lambda _{0} \in \Lambda) 
(\forall \lambda \in \Lambda)\
\lambda \geq \lambda_{0}&
\Rightarrow&
f_{\lambda}^{-1}(V^{\rm c}) \cap U = \emptyset\\
& 
\Rightarrow& 
V^{\rm c} \cap f_{\lambda }(U) = \emptyset.
\end{array}$$
Thus 
$$(\exists U \in \mathcal{N}(x)) 
(\exists \lambda_{0}\in \Lambda)
(\forall \lambda \in \Lambda) \
\lambda \geq \lambda _{0} \Rightarrow f_{\lambda}(U)\subset V.$$
By Theorem 2, $f_{\lambda}\buildrel cc\over {\longrightarrow f}$.
\end{proof}

\begin{remark}
In Theorem 3 the condition ``for every $B\subset Y$'' can be replaced by: 
for every $B=V^{\rm c}$, where $V$ is a neighbourhood basic element.
\end{remark}

\begin{corollary}
Let $(f_{\lambda })_{\lambda \in \Lambda }$ be a net in $Y^{X}$ and $Y$ be
a topological space. 
The net $(f_{\lambda})_{\lambda \in \Lambda} $ converges continuously to
$f\in Y^{X}$ if and only if the following holds:
\begin{itemize}
\item[(4*)]
$\overline{\lim\limits_{\Lambda}} f_{\lambda}^{-1}(B) \subset f^{-1}(B)$ 
for every closed subset $B$ in $Y$.
\end{itemize}
\end{corollary}

\begin{proof}
Let a net $(f_{\lambda})$ converges continuously to $f\in Y^{X}$ and
$B$ be a closed set in $Y$. 
By (4), 
$$\overline{\lim \limits_{\Lambda}} f_{\lambda}^{-1}(B) \subset
f^{-1}(\overline B) =
f^{-1}(B).$$
Conversely, suppose that (4*) holds and let $B$ be a subset in $Y$. 
By (4*), for the closed subset $\overline B$, 
$\overline{\lim\limits_{\Lambda}} f_{\lambda}^{-1}(\overline B) \subset
f^{-1}(\overline B)$, 
and by isotony of the upper limit, 
$\overline{\lim\limits_{\Lambda}} f_{\lambda}^{-1}(B) \subset
\overline {\lim\limits_{\Lambda}} f_{\lambda}^{-1}(\overline B)$. 
Hence the statement.
\end{proof}

\begin{theorem}
Let $\sigma $ and $\sigma'$ be two closure operators on $Y^{X}$. 
\begin{enumerate}
\item
If $\sigma'$ is finer than $\sigma$ and $\sigma'$ is proper, then 
$\sigma$ is proper. 
\item 
If $\sigma \leq \sigma'$ and $\sigma$ is admissible, then $\sigma'$ is 
admissible. 
\item 
If $\sigma$ is proper and $\sigma'$ is admissible, then 
$\sigma \leq \sigma'$. 
\item 
If there is a closure operator on $Y^{X}$ which is both proper and
admissible, it is unique.
\end{enumerate}
\end{theorem}

\begin{proof}
(3).
If $\sigma'$ is admissible, the evaluation mapping 
$$\varepsilon:(Y^{X},\sigma')\times (X,u) \rightarrow (Y,v)$$
is continuous by Theorem 1. 
Since $\sigma$ is proper, for $Z = ({Y}^{X},\sigma')$, the identity
$$\varepsilon^{\ast}=1_{Y^{X}}:(Y^{X},\sigma') \rightarrow 
(Y^{X},\sigma)$$
is continuous, hence $\sigma'$ is finer than $\sigma$.
\end{proof}

\begin{corollary} 
Let $\{\sigma_{\alpha}\}$ be a collection of closure operators on $Y^{X}$. 
\begin{enumerate}
\item
If $\sigma_{\alpha }$ is proper for every $\alpha$, then the infimum and
supremum, $\land \sigma_{\alpha}$ and $\lor \sigma_{\alpha}$, are
proper.
\item
If $\sigma_{\alpha}$ is admissible for every $\alpha$, then the supremum,
$\lor \sigma_{\alpha}$ is admissible as well.
\end{enumerate}
\end{corollary}

\begin{proof}
For the nontrivial part of (1), let $\sigma_{\alpha }$ be proper 
for every $\alpha$. 
That is, for any space $(Z,w)$, continuity of 
$g:(Z,w)\times (X,u) \rightarrow (Y,v)$ implies 
$g^{\ast}:(Z,w) \rightarrow (Y^{X},\sigma_{\alpha})$ 
is continuous. 
In order to prove continuity of 
$g^{\ast}:(Z,w) \rightarrow (Y^{X}, \lor\sigma_{\alpha})$,
for any $z\in Z$ and every neighbourhood $G$ of $g^{\ast}(z)$, there are
finitely many $G_{i} \in \mathcal{N}_{\alpha_{i}}(g^{\ast}(z))$ such that
$\bigcap G_{i}\subset G$. 
By continuity of 
$g^{\ast}:(Z,w) \rightarrow (Y^{X},\sigma_{\alpha})$, 
for every $i$ there is a $W_{i}\in \mathcal{N}(z)$ such that 
$g^{\ast}(W_{i})\subset G_{i}$. 
Then 
$$W = \bigcap W_{i} \in \mathcal{N}(z)$$ 
and 
$$g^{\ast}(W) \subset \bigcap g^{\ast}(W_{i}) \subset \bigcap G_{i} 
\subset G$$
holds. 
Thus $g^{\ast}:(Z,w) \rightarrow (Y^{X}, \lor\sigma_{\alpha})$ is
continuous at $z$.
\end{proof}

\begin{theorem}
A closure operator $\sigma$ on $Y^{X}$ is proper (splitting) if and only 
if continuous convergence of a net implies its convergence in 
$(Y^{X}, \sigma)$, and it is admissible (jointly continuous) if and only 
if the reverse holds, that is, convergence of a net in $(Y^{X},\sigma)$
implies its continuous convergence.
\end{theorem}

\begin{theorem}
A closure operator $\sigma $ on $Y^{X}$ is: 
\begin{enumerate}
\item
proper if and only if continuity of a mapping 
$g:Z\times (X,u) \rightarrow (Y,v)$ implies continuity of the mapping
$g^{\ast}:Z \rightarrow (Y^{X},\sigma)$ for every topological space $Z$
being either a ${\rm T}_1$-space having at most one non-isolated point or
the Sierpinski space;
\item
admissible if and only if the reverse holds, that is, continuity of a
mapping $g^{\ast}:Z \rightarrow (Y^{X},\sigma)$ implies continuity of the 
mapping $g:Z\times (X,u) \rightarrow (Y,v)$ for every topological space 
$Z$ being either a ${\rm T}_1$-space having at most one non-isolated point
or the Sierpinski space.
\end{enumerate}
\end{theorem}

\begin{corollary}
In Theorem 6 the conditions on the space $Z$ can be replaced by: 
$Z$ is a topological space having at most one non-isolated point.
\end{corollary}

In the sequel the finest proper topology on $Y^{X}$, which exists by
Corollary 2, is characterized by means of convergence classes and upper
limits, analogously to the topological situation. (Cf. \cite{IP}.) 

Let $\mathcal{C}(\sigma)$ be the {\it convergence class} of the closure
space $(Y^{X},\sigma)$, that is 
$$\mathcal{C}(\sigma) = 
\{ ((f_{\lambda})_{\lambda \in \Lambda}, f) \mid 
f_{\lambda}, f \in Y^{X} \mbox{ and } 
f_{\lambda}\buildrel \sigma\over {\longrightarrow f}
\}.$$ 
It can be easily seen that $\mathcal{C}(\sigma)$ satisfies the following
axioms: (cf.\ \cite[35 A.2. Theorem]{vC} and \cite[2.9 Theorem]{Ke})
\begin{description}
\item[(i) (CONSTANTS)]
If $(f_{\lambda })_{\lambda \in \Lambda}$ is a net such that 
$f_{\lambda} = f$ for every $\lambda \in \Lambda$, then
$(f_{\lambda})$ converges to $f$, that is, 
$((f_{\lambda})_{\lambda \in \Lambda}, f) \in \mathcal{C}(\sigma)$;
\item[(ii) (SUBNETS)]
If a net $(f_{\lambda })$ converges to $f$, so does each subnet of
$(f_{\lambda})$, i.e.\ if 
$((f_{\lambda})_{\lambda \in \Lambda}, f) \in \mathcal{C}(\sigma)$, 
then $((g_{\mu})_{\mu \in {\rm M}}, f) \in \mathcal{C}(\sigma)$ for every
subnet $(g_{\mu})$ of $(f_{\lambda})$;
\item[(iii) (DIVERGENCE)]
If a net $(f_{\lambda })$ does not converge to $f$, then there is a subnet
$(g_{\mu})$ of $(f_{\lambda })$ no subnet of which converges to $f$, 
i.e.\ 
$((f_{\lambda })_{\lambda \in \Lambda},f) \notin \mathcal{C}(\sigma)$, 
then there is a subnet $(g_{\mu})$ of $(f_{\lambda})$ such that 
$((h_{\nu})_{\nu \in {\rm N}},f) \notin \mathcal{C}(\sigma)$ 
for every $(h_{\nu})$ subnet of $(g_{\mu})$.
\end{description}

\begin{proof}
(iii).
Let $(f_{\lambda})$ be a net in $Y^{X}$, $f\in Y^{X}$ such that 
$((f_{\lambda})_{\lambda \in \Lambda},f)\notin \mathcal{C}(\sigma)$. 
It means that 
$$(\exists G_0 \in \mathcal{N}(f)) 
(\forall \lambda \in \Lambda) 
(\exists \lambda' \in \Lambda) \lambda' 
\geq 
\lambda \land f_{\lambda'} \notin G_{0}.$$
Thus there is a cofinal subset ${\rm M} \subset \Lambda $ such that
$f_{\mu }\notin G_{0}$ for every $\mu \in {\rm M} $. 
$(f_{\mu })_{\mu \in {\rm M} }$ is a subnet of $(f_{\lambda })$, no
subnet of which converges to $f$.
\end{proof}

The space $(Y^{X}, \sigma)$ is topological if and only if its convergence
class satisfies the axiom of (ITERATED LIMITS). 
(Cf.\ \cite[2.9 Theorem]{Ke} and \cite[15 B.13. and 35 A.3. Theorems]{vC}.)

Denote by $\mathcal{C}^{\ast}$ the class of all pairs 
$((f_{\lambda})_{\lambda \in \Lambda}, f)$ such that $(f_{\lambda})$ is a 
net in $Y^{X}$ which converges continuously to $f\in Y^{X}$,
i.e.\ 
$$\mathcal{C}^{\ast} = 
\{ ((f_{\lambda})_{\lambda \in \Lambda}, f) \mid 
f_{\lambda}, f \in Y^{X} \mbox{ and }
f_{\lambda}\buildrel cc\over {\longrightarrow f} 
\}.$$
By Theorem 5, $\sigma$ is proper if and only if 
$\mathcal{C}^{\ast} \subset \mathcal{C}(\sigma)$ 
and it is admissible if and only if the reverse inclusion holds: 
$\mathcal{C}(\sigma) \subset \mathcal{C}^{\ast}$.

\begin{theorem}
The class $\mathcal{C}^{\ast}$ satisfies the axioms (CONSTANTS), 
(SUBNETS) and (DIVERGENCE). 
\end{theorem}

\begin{proof}
For (CONSTANTS) and (SUBNETS) is clear. 

(DIVERGENCE).
Let $(f_{\lambda})$ be a net in $Y^{X}$, $f\in Y^{X}$ and let 
$$((f_{\lambda})_{\lambda \in \Lambda}, f) \notin \mathcal{C}^{\ast}.$$
By Theorem 3 there is a subset $B\subset Y$ such that 
$\overline{\lim\limits_{\Lambda}} f_{\lambda}^{-1}(B)$ is not contained
in $f^{-1}(v(B))$. 
Let 
$$x \in \overline{\lim\limits_{\Lambda}} f_{\lambda}^{-1}(B) \setminus 
f^{-1}(v(B)).$$
Let $\mathcal{N}(x)$ be the set of all neighbourhoods of $x$ directed by
inverse inclusion and let ${\rm M} = \Lambda \times \mathcal{N}(x)$. 
If $\mu = (\lambda ,U) \in \Lambda \times \mathcal{N}(x)$, let 
$\varphi:{\rm M} \rightarrow \Lambda$ be defined by 
$\varphi(\mu) \in \Lambda$ such that 
$\varphi(\mu) = \varphi(\lambda, U) \geq \lambda$ and 
$f_{\varphi(\mu)}^{-1}(B) \cap U \ne \emptyset$.
The net $(g_{\mu})_{\mu \in {\rm M}}$, where 
$g_{\mu} = f_{\varphi(\mu)}$, is a subnet of $(f_{\lambda})$.

Let $(h_{\nu})$ be a subnet of $(g_{\mu})$ and 
$\psi :{\rm N} \rightarrow {\rm M} $ be the corresponding map. 
In order to prove that 
$((h_{\nu})_{\nu \in {\rm N}}, f) \notin \mathcal{C}^{\ast}$, 
let $\nu_0 \in {\rm N} $ and $U \in \mathcal{N}(x)$. 
If $\psi(\nu_0) = (\lambda_0,U_0) \in {\rm M}$, set 
$\hat U = U_{0} \cap U \in \mathcal{N}(x)$ and 
$\mu_0 = (\lambda_0,\hat U)$. 
There is ${\nu_1\in {\rm N}}$ such that 
$\nu_{1}\geq \nu_0$ and 
$\nu \geq \nu_{1} \Rightarrow \psi (\nu) \geq \mu_0$. 
For any $\nu \geq \nu _{1}$ and $\psi (\nu) = (\lambda, \tilde U)$ we have
$$h_{\nu }^{-1}(B) \cap U = 
f_{\varphi(\psi(\nu))}^{-1}(B) \cap U \supset
f_{\varphi(\psi(\nu))}^{-1}(B) \cap \hat U \supset 
f_{\varphi(\psi (\nu))}^{-1}(B) \cap \tilde U \ne 
\emptyset.$$ 
It means that 
$x \in \overline{\lim\limits_{{\rm N}}} h_{\nu}^{-1}(B)$ and hence
$\overline{\lim\limits_{{\rm N}}} h_{\nu }^{-1}(B)$ is not contained in 
$f^{-1}(v(B))$. 
Thus the axiom (DIVERGENCE) is satisfied.
\end{proof}

\begin{corollary}
$\mathcal{C}^{\ast }$ is the convergence class of the finest proper
topology on $Y^{X}$ if and only if $\mathcal{C}^{\ast}$ satisfies the
axiom (ITERATED LIMITS).
\end{corollary}

\begin{theorem}
A subset $G$ in $Y^{X}$ is open in the finest proper topology if and only
if for every $f\in G$ and for every net 
$(f_{\lambda})_{\lambda \in \Lambda}$ in $Y^{X}$ such that {\rm (4)}
holds, there exists a $\lambda _{0} \in \Lambda $ such that
$f_{\lambda}\in G$ for every $\lambda {\geq \lambda}_{0}$.
\end{theorem}

\begin{proof}
($\Leftarrow $:) 
Let $\tau$ be the collection of subsets in $Y^{X}$ with the given
property.

$\tau$ is a topology on $Y^{X}$:
for $G_{1}, G_{2} \in \tau$, $f \in G_{1} \cap G_{2}$ and a net
$(f_{\lambda})_{\lambda \in \Lambda}$ satisfying (4), there are 
$\lambda_{1}, \lambda_2 \in \Lambda$ such that 
$f_{\lambda} \in G_i$ for all $\lambda \geq \lambda_{i}$, $i=1,2$. 
Then $f_{\lambda} \in G_{1} \cap G_{2}$ for all 
$\lambda \geq \lambda_{0} = \max \limits \{ \lambda_{1}, \lambda_{2} \}$. 
It follows that $G_{1} \cap G_{2} \in \tau$. 
Similarly, if $\{ G_{\alpha} \} \subset \tau$ and 
$G = \bigcup_{\alpha} \{G_{\alpha}\}$, let $f\in G$ and 
$(f_{\lambda})_{\lambda \in \Lambda}$ be a net which satisfies (4). 
There is an $\alpha _{0}$ such that $f\in G_{\alpha _{0}}$, and since
(4) holds, there exists a ${\lambda_0 \in \Lambda}$ such that 
$f_{\lambda}\in G_{\alpha_{0}}\subset G$ for every 
$\lambda \geq \lambda_{0}$. 
Thus $G \in \tau$.

$\tau$ is proper: 
let $(f_{\lambda })_{\lambda \in \Lambda}$ be a net such that 
$f_{\lambda}\buildrel cc\over {\longrightarrow f}$. 
By Theorem 3, 
$\overline{\lim\limits _{\Lambda }} f_{\lambda}^{-1}(B) \subset
f^{-1}(v(B))$ 
for every subset $B$ in $Y$. 
Let $f\in G\in \tau$. 
By the assumption, there exists a ${\lambda _0 \in \Lambda}$ such that
$f_{\lambda} \in G$ for every $\lambda \geq \lambda_{0}$, that is 
$f_{\lambda}\buildrel \tau\over {\longrightarrow f}$. 
By Theorem 5, $\tau$ is proper.

$\tau$ is the finest proper topology: 
let $\sigma$ be a proper topology on $Y^{X}$ and $H \in \sigma$.
Let $f\in H$ and a net $(f_{\lambda })_{\lambda \in \Lambda}$ satisfy (4). 
By Theorem 3, 
$f_{\lambda}\buildrel cc\over {\longrightarrow f}$. 
Since $\sigma$ is proper, 
$f_{\lambda }\buildrel \sigma\over {\rightarrow f}$. 
By definition of convergence, there exists a ${\lambda_0 \in \Lambda}$ 
such that $f_{\lambda} \in H$ for every $\lambda \geq \lambda_{0}$. 
By definition of $\tau$, $H \in \tau$. 
Thus $\sigma \subset \tau$.

($\Rightarrow $:) 
Let a subset $G$ in $Y^{X}$ be open in the finest proper topology 
$\tau$, let $f\in G$ and $(f_{\lambda })$ be a net satisfying (4). 
By Theorem 3, 
$f_{\lambda }\buildrel cc\over {\longrightarrow f}$. 
Since $\tau$ is proper, 
$f_{\lambda} \buildrel \tau\over {\longrightarrow f}$. 
By definition of convergence, there exists a ${\lambda_0 \in \Lambda}$
such that $f_{\lambda }\in G$ for every $\lambda \geq \lambda_{0}$.
\end{proof}

In order to give nontrivial examples of admissible and proper topologies
and to get results analogous to Theorems 4.1 and 4.21 in \cite{AD}, we
consider the following sets. 
Let 
$$\mathcal{V} = \{ V\subset Y | {\rm int}_{v}V \ne \emptyset \}.$$ 
For $A \subset X$ and $V \in \mathcal{V}$, set 
$$(A,V) = \{ f\in Y^{X} | f(A)\subset V \}.$$

Let $\mathcal{A}$ be a family of subsets of $X$. 
The collection 
$$
\{ (A,V) | A \in \mathcal{A}, V \in \mathcal{V}, V = {\rm int}_{v}V \},
$$
is a subbase for a topology on $Y^{X}$, which will be called the
$\mathcal{A}$-{\it topology}.

Let $\mathcal{C}$ be an interior cover of $X$. 
The collection 
$$\{ (u(K),V) | V \in \mathcal{V} \mbox{ and } K \subset X
\mbox{ is such that } u(K) \subset C 
\mbox{ for some } C \in \mathcal{C} \},$$
is a subbase for a topology on $Y^{X}$, which will be called the 
$\mathcal{C}$-{\it topology}.

\begin{theorem}
Let $(X,u)$ be a regular closure space and $(Y,v)$ be arbitrary. 
For every interior cover $\mathcal{C}$ of $X$, the $\mathcal{C}$-topology
is admissible.
\end{theorem}

\begin{proof}
By Theorem 1, it is enough to prove that the evaluation mapping is
continuous.
Let $f \in Y^{X}$, $x \in X$ and $V \in \mathcal{N}(f(x))$,
$f(x) = \varepsilon(f,x)$. 
By continuity of $f$, the set $U = f^{-1}(V) \in \mathcal{N}(x)$. 
Choose a $C \in \mathcal{C}$ so that $x \in {\rm int}_{u}C$. 
Then $U \cap C \in \mathcal{N}(x)$ and by regularity of $X$, there is a 
$U_{1} \in \mathcal{N}(x)$ such that 
$$x \in {\rm int}_{u}U_{1} \subset 
U_{1} \subset 
u(U_{1}) \subset U \cap C.$$
For the subbasic element $(u(U_1),V)$ in the $\mathcal{C}$-topology,
$\varepsilon((u(U_1),V),U_1) \subset V$ since for every 
$f_{1} \in (u(U_{1}),V)$ and each $x \in U_1$, 
$\varepsilon(f_1,x_1) = f_1(x_1) \in V$.
\end{proof}

\begin{theorem}
Let $(X,u)$ and $(Y,v)$ be closure spaces and $\mathcal{A}$ be a 
collection of compact subsets in $(X,u)$. 
The $\mathcal{A}$-topology is always proper.
\end{theorem}

\begin{proof}
Let $g:(Z,w)\times (X,u) \rightarrow (Y,v)$ be a continuous function. 
In order to prove continuity of the mapping 
$g^{\ast}:(Z,w) \rightarrow (Y^{X},\sigma)$, where $Y^{X}$ is endowed
with the $\mathcal{A}$-topology, let $z\in Z$ and $f=g^{\ast}(z)$. 
For a subbasic element $(K,V)$ containing $f$, where $K$ is a compact set
in $X$, 
$$f(K) = g(\{z\} \times K) \subset {\rm int}_{v}V = V.$$
By continuity of $g$, 
$$(\forall x \in K) \ 
g(z,x) \in {\rm int}_{v}V 
\Leftrightarrow 
V \in \mathcal{N}(g(z,x))$$
implies
$$(\forall x\in K) 
(\exists W_{x}\in \mathcal{N}(z)) 
(\exists U_{x}\in \mathcal{N}(x)) \ g(W_{x}\times U_{x}) \subset V.$$

$(\forall x \in K) U_x \in \mathcal{N}_{x}$
implies 
$\{ U_{x} \mid x \in K\}$ 
is an interior cover of the compact set $K$, so there is a finite subcover
$\{U_{x_{i}} \mid i=1,\cdots k\}$. 
Set $W = \bigcap_{i=1}^{k}W_{x_{i}}$. 
Then $W \in \mathcal{N}(z)$. 
It follows that 
$$g(W \times K) \subset 
g(\bigcup_{i=1}^{k}{(W}_{x_{i}}\times U_{x_{i}})) \subset 
V.$$ 
Thus 
$$(\forall z' \in W) \
g(z',K) \subset V \Rightarrow g^{\ast}(W) \subset (K,V).$$
\end{proof}

\section{Some special cases including $\theta$-closure}

\subsection{}
Let $(X,u)$ be a closure space and $(Y,\mathcal{V})$ be a topological
space. 
A function $f:(X,u) \rightarrow (Y,\mathcal{V})$ is continuous if and
only if the function $f:(X,\hat u) \rightarrow (Y, \mathcal{V})$
between topological spaces is continuous, where $\hat u$ is the
topological modification of the closure operator $u$. 
In that case the problem is reduced to the topological case since 
$Y^{X} = C((X,\hat u),Y)$ 
and the topological modification of the product of closure spaces is the
product of topological modifications.

\subsection{}
It was already remarked that $\theta$-continuous functions are continuous
functions of the corresponding \v Cech closure spaces. 
Compact sets in $(X,{\rm cl}_{\theta })$ are (quasi-)H-closed 
(g-H-closed) in $(X,\mathcal{U})$. 
Thus Theorems 3.1--3.6 and 4.2 in \cite{C} are special cases of Theorems
2, 1, 4, 5 and 10 respectively.

If $(Y,\mathcal{V})$ is regular, the $\theta$-topology 
$\mathcal{V}_{\theta} = \mathcal{V}$. 
Then for a topological space $(X,\mathcal{U})$, a function
$f:(X,\mathcal{U}) \rightarrow (Y,\mathcal{V})$ is 
$\theta$-continuous if and only if 
$f:(X,{\rm cl}_{\theta}) \rightarrow (Y,\mathcal{V})$ 
is continuous, which is equivalent to 
$f:(X,\mathcal{U}_{\theta}) \rightarrow (Y,\mathcal{V})$ 
be continuous \cite[Thm 16.B.4]{vC}. 
Note that the topological modification of 
$\mathcal{U}{\rm cl}_{\theta }$ is $\mathcal{U}_{\theta}$, 
the topology of $\theta$-open sets in $(X,\mathcal{U})$.

\subsection{}
A large number of continuous-like mappings between topological spaces is 
known in the literature. 
Recently, Georgiou and Papadopoulos \cite{GP1, GP2} considered some
of them and investigated splitting and jointly continuous topologies on
the sets of these functions. 
Let us remark that all these examples and the main results are special
cases of our subjects of investigations. 
For, let $(X,\mathcal{U})$ and $(Y,\mathcal{V})$ be topological spaces. 
It follows from the definitions that a function 
$f:(X,\mathcal{U})\rightarrow (Y,\mathcal{V})$ is:
\begin{itemize}
\item[(1)]
{\it strongly $\theta $-continuous} (cf.\ \cite{GP1}) at a point $x$
(on the
set $X$) if and only if 
$f:(X,\mathcal{U}{\rm cl}_{\theta}) \rightarrow (Y,\mathcal{V})$ 
is continuous at $x$ (on the set $X$);
\item[(2)]
{\it super-continuous} (cf.\ \cite{GP2}) at $x$ (on the set $X$) if
and only if 
$f:(X,\mathcal{U}_s) \rightarrow (Y,\mathcal{V})$ 
is continuous at $x$ (on the set $X$), where $\mathcal{U}_s$ is the 
semi-regularization topology of $\mathcal{U}$ (see \cite{GMRV} for
example);
\item[(3)]
{\it weakly continuous} (cf.\ \cite{GP2}) at $x$ (on the set $X$) if
and only if 
$f:(X,\mathcal{U}) \rightarrow (Y,\mathcal{V}{\rm cl}_{\theta})$
is continuous at $x$ (on the set $X$);
\item[(4)]
{\it weakly $\theta $-continuous} (cf.\ \cite{GP2}) at a point $x$ (on
$X$) if and only if 
$f:(X,\mathcal{U}_s) \rightarrow (Y,\mathcal{V}{\rm cl}_{\theta})$ 
is continuous at $x$ (on $X$).
\end{itemize}

Also $\theta$-{\it convergence} of a net $(x_{\mu})$ in $(X,\mathcal{U})$
means convergence of $(x_{\mu})$ in the corresponding closure space
$(X,\mathcal{U}{\rm cl}_{\theta})$, while {\it weak $\theta$-convergence} 
of a net $(x_{\mu})$ in $(X,\mathcal{U})$ is convergence of $(x_{\mu})$ in
$(X,\mathcal{U}_s)$.

Similarly, {\it $\theta $-continuous convergence} (respectively: 
{\it strongly $\theta$-continuous convergence, weakly $\theta$-continuous 
convergence, weakly continuous convergence} and {\it super continuous
convergence}) of a net $(f_{\lambda })$ in $Y^X$ is continuous convergence
of $(f_{\lambda })$ for the corresponding closure spaces. 
Thus we are concerned with a change of topology, better to say: change of
closure operator, technique. 
So the main results in \cite{GP1} and \cite{GP2} are special cases of
the above Theorems 1--10.

%\bibliographystyle{amsplain}
%\bibliography{21}
\providecommand{\bysame}{\leavevmode\hbox to3em{\hrulefill}\thinspace}
\providecommand{\MR}{\relax\ifhmode\unskip\space\fi MR }
\providecommand{\MRhref}[2]{%
  \href{http://www.ams.org/mathscinet-getitem?mr=#1}{#2}
}
\providecommand{\href}[2]{#2}

\end{document}